\documentclass[12pt]{article}
\input epsf.tex

\usepackage{amsmath}
\usepackage{amsthm}
\usepackage{amsfonts}
\usepackage{amssymb}
\usepackage{graphicx}
\usepackage{latexsym}

\usepackage{amsmath,amsthm,amsfonts,amssymb,mathrsfs}

\usepackage{epsfig}

\theoremstyle{plain}
\newtheorem{thm}{Theorem}[section]
\newtheorem{prop}[thm]{Proposition}
\newtheorem{lem}[thm]{Lemma}
\newtheorem{cor}[thm]{Corollary}

\theoremstyle{definition}
\newtheorem{defn}{Definition}
\theoremstyle{remark}

\advance \topmargin by -\headheight
\advance \topmargin by -\headsep
\evensidemargin \oddsidemargin
\def\cc{{\curvearrowright}}

\def\F{{\mathbb F}}

\def\cG{{\mathcal G}}

\def\rng{{\textrm{rng }}}

\def\chix{{\raise.5ex\hbox{$\chi$}}}

\def\Z{{\mathbb Z}}

\begin{document}
\title{Orbit equivalence, coinduced actions and free products}
\author{Lewis Bowen}
\begin{abstract}
The following result is proven. Let $G_1 \cc^{T_1} (X_1,\mu_1)$ and $G_2 \cc^{T_2} (X_2,\mu_2)$ be orbit equivalent (OE), essentially free, probability measure preserving actions of countable groups $G_1$ and $G_2$.  Let $H$ be any countable group. For $i=1,2$, let $\Gamma_i = G_i *H$ be the free product. Then the actions of $\Gamma_1$ and $\Gamma_2$ coinduced from $T_1$ and $T_2$ are OE. As an application, it is shown that if $\Gamma$ is a free group, then all nontrivial Bernoulli shifts over $\Gamma$ are OE.
\end{abstract}
\maketitle
\noindent
{\bf Keywords}: coinduced, orbit equivalence, Bernoulli shifts.\\
{\bf MSC}:37A20\\

\noindent


\section{Introduction}

Let $G$ be a countable group and $(X,\mu)$ a standard probability space. A probability measure-preserving (p.m.p.) action of $G$ on $(X,\mu)$ is a collection $\{T^g\}_{g\in G}$ of measure-preserving transformations $T^g:X \to X$ such that $T^{g_1}T^{g_2}=T^{g_1g_2}$ for all $g_1,g_2 \in G$. We denote this by $G \curvearrowright^T (X,\mu)$.

Suppose $G_1\cc^{T_1} (X_1,\mu_1)$ and $G_2\cc^{T_2} (X_2,\mu_2)$ are two  p.m.p. actions. A measurable map $\phi:X_1' \to X_2'$ (where $X'_i \subset X_i$ is conull) is an {\em orbit equivalence} if the push-forward measure $\phi_*\mu_1$ equals $\mu_2$ and for every $x\in X_1'$, $\{T_1^gx:~g\in G_1\}=\{T_2^g\phi(x):~g\in G_2\}$. If there exists such a map, then the actions $T_1$ and $T_2$ are said to be {\em orbit equivalent} (OE). If, in addition, there is a group isomorphism $\Psi:G_1 \to G_2$ such that $\phi(T_1^g x) = T_2^{\Psi(g)} \phi(x)$ for every $x\in X'_1$ and $g\in G_1$ then the actions $T_1$ and $T_2$ are said to be {\em measurably-conjugate}.



The initial motivation for orbit equivalence comes from the study of von Neumann algebras. It is known that two p.m.p. actions are orbit equivalent if and only if their associated crossed product von Neumann algebras are isomorphic by an isomorphism that preserves the Cartan subalgebras [Si55].  H. Dye [Dy59, Dy63]  proved the pioneering result that any two ergodic p.m.p. actions of the group of integers on the unit interval are OE. This was extended to amenable groups in [OW80] and [CFW81]. By contrast, it is now known that every nonamenable group admits a continuum of non-orbit equivalent ergodic p.m.p. actions [Ep09]. This followed a series of earlier results that dealt with various important classes of non-amenable groups ([GP05], [Hj05], [Io09], [Ioxx], [Ki08], [MS06], [Po06]). 

In the last decade, a number of striking OE rigidity results have been proven (for surveys, see [Fu09], [Po07] and [Sh05]). These imply that, under special conditions, OE implies measure-conjugacy. By contrast, the main theorem of this paper is a contribution to the relatively small number of OE ``flexibility'' results. 

To illustrate the new results, let us consider the classification of Bernoulli shifts up to orbit equivalence and measure-conjugacy. So let $G$ be a countable group. Let $(K,\kappa)$ be a standard probability space. $K^G$ is the set of all of functions $x:G \to K$ with the product Borel structure. For each $g\in G$, let $S^g:K^G \to K^G$ be the shift-map defined by $S^gx(h):=x(g^{-1}h)$ for any $h\in G$ and $x\in K^G$. This map preserves the product measure $\kappa^G$. The action $G \cc^S (K^G,\kappa^G)$ is called the {\em Bernoulli shift over $G$ with base-space $(K,\kappa)$}.

If $\kappa$ is supported on a finite or countable set $K'\subset K$ then the {\em entropy} of $(K,\kappa)$ is defined by
$$H(K,\kappa) := -\sum_{k \in K'} \kappa(\{k\}) \log\big( \kappa(\{k\})\big).$$
Otherwise, $H(K,\kappa):=+\infty$.

A. N. Kolmogorov proved that if two Bernoulli shifts $\Z \cc (K^\Z,\kappa^\Z)$ and $\Z \cc (L^\Z,\lambda^\Z)$ are measurably-conjugate then the base-space entropies $H(K,\kappa)$ and $H(L,\lambda)$ are equal [Ko58, Ko59]. This answered a question of von Neumann which had been posed at least 20 years prior. The converse to Kolmogorov's theorem was famously proven by D. Ornstein [Or70ab]. Both results were extended to countable infinite amenable groups in [OW87].  

A group $G$ is said to be {\em Ornstein} if whenever $(K,\kappa), (L,\lambda)$ are standard probability spaces with $H(K,\kappa)=H(L,\lambda)$ then the corresponding Bernoulli shifts $G \cc (K^G,\kappa^G)$ and $G \cc (L^G,\lambda^G)$ are measurably conjugate. A. M. Stepin proved that if $G$ contains an Ornstein subgroup, then $G$ is Ornstein [St75]. Therefore, any group $G$ that contains an infinite amenable subgroup is Ornstein. It is unknown whether every countably infinite group is Ornstein.

By [Bo10], every sofic group satisfies a Kolmogorov-type theorem. Precisely, if $G$ is sofic, $(K,\kappa), (L,\lambda)$ are standard probability spaces with $H(K,\kappa) + H(L,\lambda)<\infty$ and the associated Bernoulli shifts $G \cc (K^G,\kappa^G)$, $G \cc (L^G,\lambda^G)$ are measurably-conjugate then $H(K,\kappa)=H(L,\lambda)$. If $G$ is also Ornstein then the finiteness condition on the entropies can be removed. Sofic groups were defined implicitly by M. Gromov [Gr99] and explicitly by B. Weiss [We00]. For example, every countably infinite linear group is sofic and Ornstein. It is unknown whether or not all countable groups are sofic. 

In summary, there is a large class of groups (e.g., all countable linear groups) for which Bernoulli shifts are completely classified up to measure-conjugacy by base-space entropy. Let us now turn to the question of orbit equivalence.

By aforementioned results of [OW80] and [CFW81], it follows that if $G_1$ and $G_2$ are any two infinite amenable groups then any two nontrivial Bernoulli shifts $G_1 \cc (K^{G_1},\kappa^{G_1})$, $G_2 \cc (L^{G_2},\lambda^{G_2})$ are OE. By contrast, the Kolmogorov-type theorem of [Bo10] combined with rigidity results of S. Popa [Po06, Po08] and Y. Kida [Ki08] prove that for many nonamenable groups $G$, Bernoulli shifts are classified up to OE by base-space entropy. For example, this includes PSL$_n(\Z)$ for $n>2$, mapping class groups of surfaces (with a few exceptions) and any nonamenable sofic Ornstein group of the form $G=H\times N$ with both $H$ and $N$ countably infinite that has no nontrivial finite normal subgroups. Until the present paper, it was an open question whether any two orbit-equivalent Bernoulli shifts over a non-amenable group $G$ are necessarily measurably-conjugate. This question featured in several talks given by S. Popa.

As a corollary to the main result we obtain:
\begin{thm}\label{thm:bernoulli}
Let $G_1, G_2$ be any two countably infinite amenable groups. Let $H$ be any countable group. For $i=1,2$, let $\Gamma_i=G_i*H$ be the free product. Let $(K,\kappa), (L,\lambda)$ be nontrivial standard probability spaces. Then the Bernoulli shifts $\Gamma_1 \cc (K^{\Gamma_1},\kappa^{\Gamma_1})$ and $\Gamma_2 \cc (L^{\Gamma_2},\lambda^{\Gamma_2})$ are orbit equivalent.
\end{thm}

It is known that amenable groups are sofic and a free product of sofic groups is sofic [ES06]. Therefore, if $H$ is sofic then each group $\Gamma_i$ above is also sofic. Since $\Gamma_i$ contains the infinite amenable group $G_i$, it is also Ornstein. Thus [Bo10] implies that the Bernoulli shifts over $\Gamma_i$ are completely classified by base-space entropy. In particular, there is a 1-parameter family of non-measurably-conjugate Bernoulli shifts over $\Gamma_i$. 

\begin{cor}
Let $G$ be a free group. Then all nontrivial Bernoulli shifts over $G$ are OE.
\end{cor}

\subsection{Statement of results}

To formulate the main theorem, we need to discuss co-induced actions (which have previously appeared for various purposes in [Lu00], [Ga05], [Da06], [DGRS08] and [Ioxx]).


\begin{defn}\label{defn:coinduced}
Fix a countable group $\Gamma$ and a subgroup $G<\Gamma$. Let $G \cc^T (X,\mu)$ be a measure-preserving action of $G$ on a standard probability space. Let $X^\Gamma$ be the space of all maps $f:\Gamma\to X$ with the product Borel structure. For each $\gamma \in \Gamma$, $S^\gamma:X^\Gamma \to X^\Gamma$ is defined by 
$$S^\gamma f(\gamma_0) := f(\gamma^{-1}\gamma_0)~ \forall f\in X^\Gamma , \gamma_0 \in \Gamma.$$
Let
$$\cG := \{f \in X^\Gamma:~ f(\gamma g) = T^{g^{-1}}f(\gamma) ~\forall \gamma \in \Gamma,g\in G\}.$$
$\cG$ is invariant under the $S$-action of $\Gamma$. We will construct a shift-invariant measure on $\cG$. To do this we need the following notion.

A {\em section} for the inclusion $G<\Gamma$ is a map $\sigma: \Gamma/G \to \Gamma$ such that $\sigma(\gamma G) \in \gamma G$ for all $\gamma\in \Gamma$. Fix such a section $\sigma$ with $\sigma(G)=e$. 

Define 
$$\Phi:\cG \to X^{\Gamma/G},~\Phi(f)(C):=f(\sigma(C))~\forall f\in\cG, C\in \Gamma/G.$$
This is a bijection. Define a measure $\nu$ on $\cG$ by $\nu(E) := \mu^{\Gamma/G}(\Phi(E))$ for all Borel $E\subset \cG$ where $\mu^{\Gamma/G}$ is the product measure on $X^{\Gamma/G}$.  An exercise reveals that $\nu$ is shift-invariant and independent of the choice of section $\sigma$. We extend $\nu$ to all of $X^\Gamma$ by setting $\nu(X^\Gamma - \cG)=0$. Then $\Gamma \cc^S (X^\Gamma,\nu)$ is called the {\em action coinduced from} $G \cc^T (X,\mu)$.
\end{defn}

\begin{thm}\label{thm:main}
Let $G_1 \cc^{T_1} (X_1,\mu_1)$ and $G_2 \cc^{T_2} (X_2,\mu_2)$ be orbit equivalent, essentially free, p.m.p. actions of countable groups $G_1$ and $G_2$.  Let $H$ be any countable group. For $i=1,2$, let $\Gamma_i = G_i *H$ be the free product. Then the co-induced actions $\Gamma_1 \cc^{S_1} (X_1^{\Gamma_1}, \nu_1)$ and  $\Gamma_2 \cc^{S_2} (X_2^{\Gamma_2}, \nu_2)$ are orbit equivalent.
\end{thm}


We can now prove theorem \ref{thm:bernoulli}.
\begin{proof}[Proof of theorem \ref{thm:bernoulli}]
Because $G_1$ and $G_2$ are countably infinite amenable groups, the well-known results of [OW80] and [CFW81] imply the Bernoulli shifts $G_1 \cc (K^{G_1},\kappa^{G_1})$ and $G_2 \cc (L^{G_2}, \lambda^{G_2})$ are orbit equivalent. Let $X=K^{G_1}$, $Y=L^{G_2}$. It is a straightforward exercise to show that for each $i=1,2$ the coinduced actions $\Gamma_1 \cc (X^{\Gamma_1},\nu_1)$ and $\Gamma_2 \cc (Y^{\Gamma_2},\nu_2)$ are measurably-conjugate to the Bernoulli shifts $\Gamma_1 \cc (K^{\Gamma_1}, \kappa^{\Gamma_1})$ and $\Gamma_2 \cc (L^{\Gamma_2},\lambda^{\Gamma_2})$. By theorem \ref{thm:main} above, these are orbit equivalent.
\end{proof}

\subsection{The idea behind the construction}
This section provides a non-rigorous sketch of the following.
\begin{thm}
Let $(X,\mu)$ be a standard probability space. Let $\Gamma_A = \langle A,C\rangle$ and $\Gamma_B=\langle B,C\rangle$ each be a free group on two generators. Let $\langle A \rangle \cc (X,\mu)$ and $\langle B \rangle \cc (X,\mu)$ be actions of the infinite cyclic groups $\langle A\rangle$ and $\langle B\rangle$ that have the same orbits. I.e., for every $x\in X$, $\{A^n x:~n\in\Z\}=\{B^n x:n\in\Z\}$. Let $\Gamma_A \cc (X^{\Gamma_A},\nu_A)$ be the action coinduced from $\langle A \rangle \cc (X,\mu)$ and let $\Gamma_B \cc (X^{\Gamma_B},\nu_B)$ be the action coinduced from $\langle B \rangle \cc (X,\mu)$. Then $\Gamma_A \cc (X^{\Gamma_A},\nu_A)$  and $\Gamma_B \cc (X^{\Gamma_B},\nu_B)$ are orbit equivalent.
\end{thm}
 This theorem is implied immediately by theorem \ref{thm:main}. The sketch we provide gives the idea behind the proof of theorem \ref{thm:main}. 
 
   

\begin{figure}[htb]
\begin{center}
\ \psfig{file=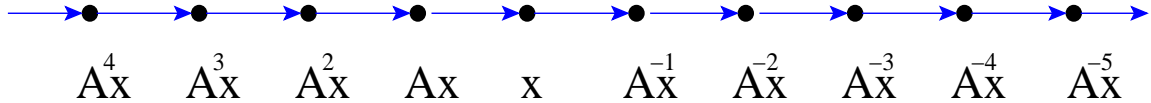,height=0.5in,width=6in}
\caption{A diagram for the $A$-orbit of $x$.}
\label{fig:1}
\end{center}
\end{figure}

Let $x\in X$ and consider the diagram of its $A$-orbit show in figure \ref{fig:1}. The vertices represent elements of the orbit. For each $i$ there is an arrow from the vertex representing $A^ix$ to the vertex representing $A^{i-1}x$. It may seem backwards to draw the arrows this way; but there is a good reason. If $X=K^{\Z}$ and $A$ is the shift map then $(A^{-n}x)(0)=x(n)$ for any $n\in\Z$. So in this case, we could replace each vertex labeled $A^{-n}x$ with the value $A^{-n}x(0)=x(n)$ and we would then have a picture of the sequence $\{x(n)\}_{n\in\Z}$. With this convention in mind, it is not necessary to label every vertex. As long as one vertex is labeled, the rest of the labels are determined by the arrows. 

\begin{figure}[htb]
\begin{center}
\ \psfig{file=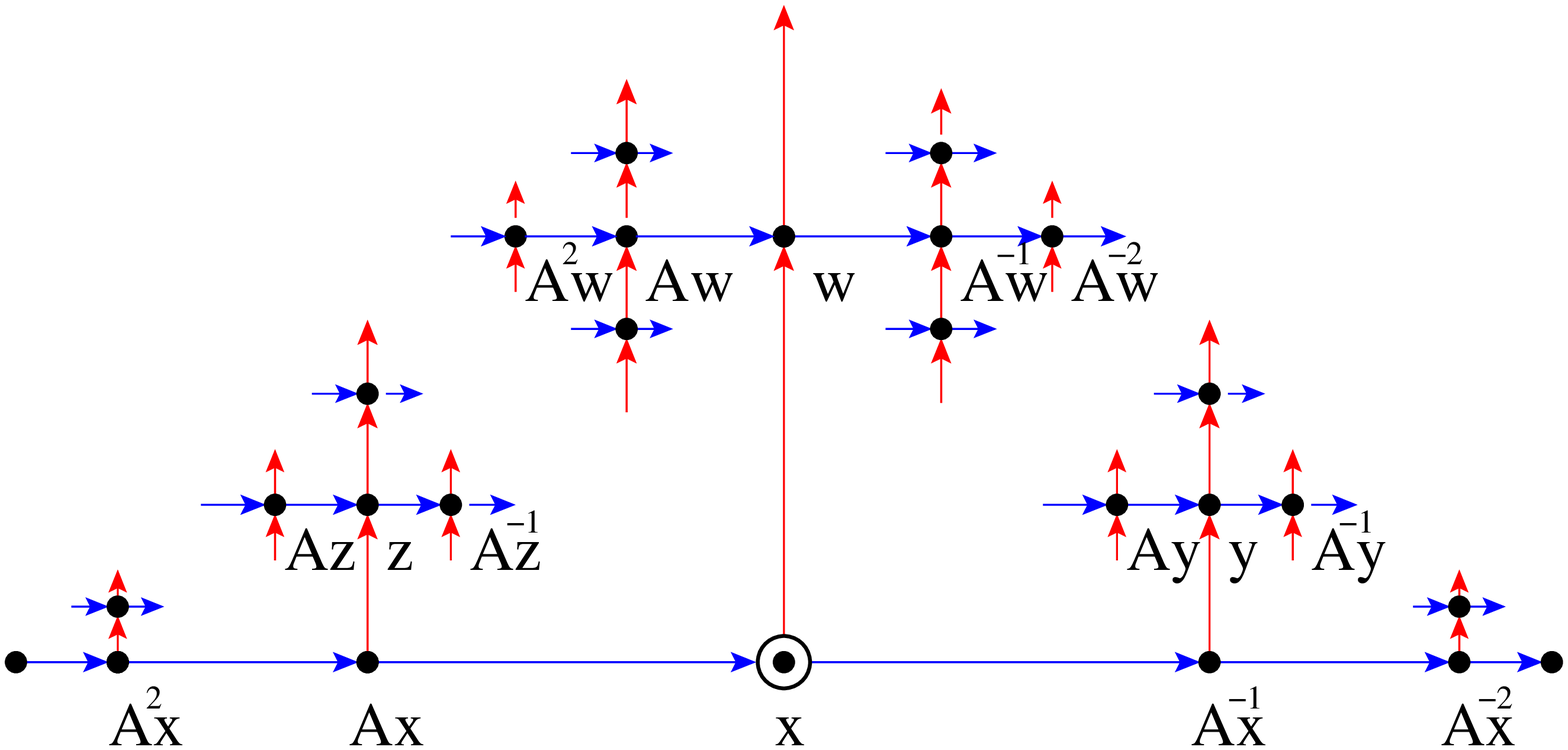,height=3in,width=6in}
\caption{A diagram for the $\langle A,C\rangle$-orbit of a function $f \in X^{\Gamma_A}$.}
\label{fig:2}
\end{center}
\end{figure}

In figure \ref{fig:2}, there is a diagram for a typical point $f \in X^{\Gamma_A}$ with respect to the coinduced measure $\nu_A$. The underlying graph is the Cayley graph of $\Gamma_A$ (only part of which is shown in the figure). The circled dot represents the identity element in $\Gamma_A$. For every $g \in \Gamma_A$ there are directed edges $(g,gA)$ and $(g,gC)$. Edges of the form $(g,gA)$ are drawn horizontally while those of the form $(g,gC)$ are drawn vertically.

Several vertices are labeled by elements $x,y,z,w \in X$. These are the values of $f$ at the group elements represented by the vertices. For examples, the diagram implies that $f(e)=x, f(AC)=y, f(A^{-1}C)=z$ and $f(C)=w$. Because $f$ is a typical point in $X^{\Gamma_A}$ (according to the measure $\nu_A$), we must have $f(\gamma A^n) = A^{-n}f(\gamma)$ for any $\gamma\in \Gamma_A$ and $n \in\Z$. This explains the other labels. Notice that if $f\in X^{\Gamma_A}$ is chosen at random with law $\nu_A$ then $x,y,z,w$ are independent samples drawn from $(X,\mu)$.



\begin{figure}[htb]
\begin{center}
\ \psfig{file=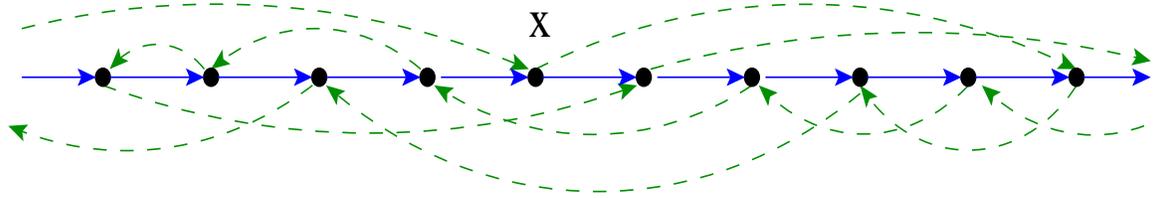,height=1in,width=6in}
\caption{The $A$ and $B$-orbit diagram of $x$.}
\label{fig:3}
\end{center}
\end{figure}

Now the $B$-orbit of a point $x\in X$ equals its $A$-orbit. So we may draw them together as in figure \ref{fig:3}. The dashed arrows represent the $B$-action. For example, one can see from the diagram that $B^{-1}x = A^{-5}x$. 

Let 
$$\cG_A := \{f \in X^{\Gamma_A}:~ f(\gamma A^n) = A^{-n}f(\gamma) ~\forall \gamma \in \Gamma_A ,n\in\Z\}$$
and define $\cG_B$ similarly. 

\begin{figure}[htb]
\begin{center}
\ \psfig{file=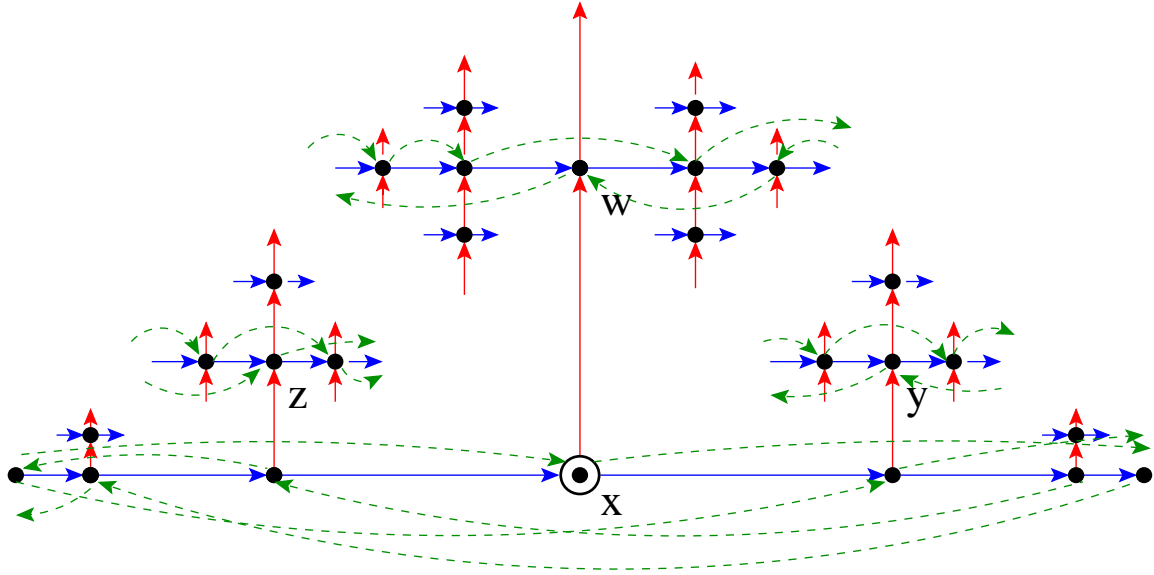,height=3in,width=6in}
\caption{A diagram for the orbit equivalence $\Omega$.}
\label{fig:4}
\end{center}
\end{figure}

The diagrams above give us the idea for how to construct the orbit equivalence between $\Gamma_A \cc (\cG_A,\nu_A)$ and $\Gamma_B \cc (\cG_B,\nu_B)$. Take the diagram of a typical point in $f\in (\cG_A,\nu_A)$ (figure \ref{fig:2}) and draw in the dashed green arrows representing the action of $\langle B\rangle$ to obtain figure \ref{fig:4}. Then erase the blue arrows and what we have left is a diagram of a point in $(\cG_B,\nu_B)$. This defines a map $\Omega:\cG_A \to \cG_B$. 

By reversing the roles of $A$ and $B$, we can similarly define a map $\Theta:\cG_B \to \cG_A$ such  that $\Omega \Theta$ and $\Theta\Omega$ are the identity maps on $\cG_B$ and $\cG_A$ respectively. So $\Omega$ is invertible. It clearly takes $\Gamma_A$-orbits to $\Gamma_B$-orbits. It might not be obvious, but $\Omega_*\nu_A = \nu_B$. Thus $\Omega$ is the required orbit equivalence. To prove theorem \ref{thm:main}, we will construct an orbit equivalence in a similar manner.

\subsection{Organization}
In \S \ref{sec:defn} we construct a map $\Omega: X^{\Gamma_1} \to X^{\Gamma_2}$. The rest of the paper is devoted to showing that this map is an orbit equivalence. There are three statements to be proven: $\Omega$ has a measurable inverse (accomplished in \S \ref{sec:invert}), $\Omega$ takes orbits to orbits (accomplished in \S \ref{sec:orbits}) and $\Omega_*\nu_1=\nu_2$, i.e., $\Omega$ is a measure-space isomorphisms. This is obtained in the last section, \S \ref{sec:measure}.

{\bf Acknowledgements}. I'd like to thank Sorin Popa for asking whether the full 2-shift over $\F_2$ is stably orbit equivalent to a Bernoulli shift over $\F_3$. My investigations into this question led to this work. I'd also like to thank Bruno Duchesne for pointing out errors in a preliminary version of this paper.

\section{Defining the orbit equivalence}\label{sec:defn}

Without loss of generality, we may assume that $X_1=X_2$, $\mu_1=\mu_2$ and the identity map from $X_1$ to $X_2$ is an orbit equivalence between $G_1 \cc^{T_1} (X_1,\mu_1)$ and $G_2 \cc^{T_2} (X_2,\mu_2)$. In other words, we may assume that $(X,\mu)$ is a standard probability space and $G_1 \cc^{T_1} (X,\mu)$, $G_2 \cc^{T_2} (X,\mu)$ are measure-preserving actions such that for a.e. $x\in X$, $\{T_1^g x:~g\in G_1\} =\{T_2^gx :~g\in G_2\}$. 

To simplify notation, we will $g\cdot x$ instead of $T_1(g)x$ if $g\in G_1, x\in X$ (or if $g\in G_2, x\in X$). Likewise we will write $\gamma f = S_\gamma f$ if $\gamma \in \Gamma_i$ and $f \in X^{\Gamma_i}$.

After reducing $X$ by removing a set of measure zero, we may assume that both actions $T_1$ and $T_2$ are free and that every $T_1$-orbit is a $T_2$-orbit and vice versa. Let $\omega:G_1 \times X \to G_2$ be the Zimmer 1-cocycle. It satisfies
\begin{eqnarray}
\omega(g,x)\cdot x &=& g\cdot x,~\forall g\in G_1, x\in X,\label{eqn:cocycle1}\\
\omega(g_1g_2,x) &=& \omega(g_1, g_2\cdot x)\omega(g_2,x),~\forall g_1,g_2\in G_1, x\in X.\label{eqn:cocycle2}
\end{eqnarray}
As in definition \ref{defn:coinduced}, for $i=1,2$ let 
$$\cG_i=\big\{f\in X^{\Gamma_i}:~ f(\gamma g) = T_i(g^{-1})f(\gamma) ~\forall \gamma \in \Gamma_i,g\in G_i\big\}.$$

\begin{prop}\label{prop:beta}
There exists a unique Borel map $\beta:\Gamma_1 \times \cG_1 \to \Gamma_2$ satisfying the following.
\begin{itemize}
\item $\beta(g,f) = \omega(g,f(e)),~ \forall g\in G_1$.
\item $\beta(h,f) = h,~\forall h\in H$.
\item $\beta(h\gamma, f) = h\beta(\gamma,f), ~\forall h\in H, \gamma \in \Gamma_1$.
\item $\beta(g\gamma,f) = \beta(g,\gamma f) \beta(\gamma,f) = \omega(g, f(\gamma^{-1}))\beta(\gamma,f), ~\forall g\in G_1,\gamma \in \Gamma_1$.

\end{itemize}
\end{prop}

\begin{proof}
Since $\gamma \in \Gamma_1 = G_1 *H$ is a free product, any element $\gamma \in \Gamma_1$ can be written uniquely as $\gamma = h_1g_1h_2g_2\cdots h_ng_n$ with $g_1, g_2,\ldots,g_{n-1} \in G_1-\{e\}$, $g_n \in G_1$, $h_1 \in H$ and $h_2,\ldots,h_{n} \in H-\{e\}$ for some $n\ge 1$. 

For $i=1\ldots n$, let $\gamma_i = h_ig_i\cdots h_ng_n$. Let $\gamma_{n+1}=e$. Define $\beta(e,f):=e$ for any $f \in X^{\Gamma_1}$. For induction, assume that $\beta(\gamma_i,f)$ has been defined for some $i\ge 2$. Then we define $\beta(\gamma_{i-1},f)$ by
$$\beta(\gamma_{i-1},f) :=h_{i-1}\omega(g_{i-1}, f(\gamma_i^{-1})) \beta(\gamma_i,f).$$
This defines $\beta$ for all $\gamma \in \Gamma_2$. It is clear that any function satisfying the conditions above must be defined as such. This explains the uniqueness part of the proposition.

It is immediate that $\beta_f$ satisfies the first three items. To check the last item, let $\gamma=h_1g_1h_2g_2\cdots h_ng_n \in \Gamma_1$ as above. Let $q \in G_1$ and $f\in X^{\Gamma_1}$. We must show that
$$\beta(q\gamma,f) = \beta(q,\gamma f)\beta(\gamma,f).$$
This is true by definition if $h_1\ne e$. So assume $h_1=e$.  Let $\tau = h_2g_2\cdots h_{n}g_{n}$ and $p=g_1$. We must show that
$$\beta(qp\tau,f) = \beta(q,p\tau  f) \beta(p\tau,f) .$$
Equivalently, we must show that 
$$\beta(qp,\tau f) \beta(\tau,f) = \beta(q,p\tau  f)\beta(p, \tau f)\beta(\tau,f).$$
Cancel the last factor on both sides and use item (1) to obtain the equivalent formulation:
$$\omega(qp, f(\tau^{-1}))  = \omega(q,f(\tau^{-1}p^{-1})) \omega(p, f(\tau^{-1})).$$
By equation (\ref{eqn:cocycle2})
$$\omega(qp,f(\tau^{-1})) = \omega(q, p\cdot f(\tau^{-1})) \omega(p,f(\tau^{-1})).$$
Since $f \in \cG_1$, $p\cdot f(\tau^{-1}) = f(\tau^{-1}p^{-1})$. This finishes the proposition.
\end{proof}

\begin{lem}\label{lem:bijective}
For every $f\in \cG_1$, $\beta(\cdot,f):\Gamma_1 \to \Gamma_2$ is a bijection.
\end{lem}
\begin{proof}
Every $\gamma \in \Gamma_2 = G_2 *H$ can be written uniquely as $\gamma = h_1g_1h_2g_2\cdots h_ng_n$ with $g_1, g_2,\ldots,g_{n-1} \in G_1-\{e\}$, $g_n \in G_1$, $h_1 \in H$ and $h_2,\ldots,h_{n} \in H-\{e\}$ for some $n\ge 1$.  For $i=1\ldots n$, let $\gamma_i = h_ig_i\cdots h_ng_n$. Let $\gamma_{n+1}=e$. 

Note that $\beta(e,f) = e=\gamma_{n+1}$. Since the action of $G_1$ on $X$ is free there exists a unique $g'_n \in G_1$ such that $\omega(g'_n,f(e)) = g_n$. Proposition \ref{prop:beta} implies  $\beta(h_n g'_n,f ) = h_ng_n=\gamma_n$. Moreover, $h_ng'_n$ is the unique element of $\Gamma_1$ with $\beta(\cdot,f)$-image equal to $\gamma_n$.

For induction, suppose that for some $i \ge 2$, $g'_j\in G_1$ has been defined for all $i \le j \le n$ and $h_ig'_i\cdots h_n g'_n$ is the unique element of $\Gamma_1$ with $\beta(\cdot,f)$-image equal to $\gamma_i$. Let $\tau =h_ig'_i\cdots h_n g'_n$. Since the action of $G_1$ on $X$ is free there exists a unique $g'_{i-1} \in G_1$ such that $\omega(g'_{i-1},f(\tau^{-1})) = g_{i-1}$. Proposition \ref{prop:beta} implies
$$\beta( h_{i-1}g_{i-1} \tau) = h_{i-1} g_{i-1} \gamma_i = \gamma_{i-1}.$$
Moreover, $h_{i-1}g_{i-1} \tau$ is the unique element of $\Gamma_1$ with $\beta(\cdot,f)$-image equal to $\gamma_{i-1}$. The completes the induction step. Since $\gamma \in \Gamma_2$ is arbitrary, this proves the lemma.
\end{proof}

Define $\Omega:\cG_1 \to \cG_2$ by 
$$(\Omega f)(\beta(\gamma^{-1},f)^{-1}) = f(\gamma), ~\forall f\in\cG_1, \gamma \in \Gamma_1.$$
 We will show that $\Omega$ is an orbit equivalence between the coinduced actions $\Gamma_1 \cc (\cG_1, \nu_1)$ and  $\Gamma_2 \cc (\cG_2, \nu_2)$. There are three statements to prove: (1) $\Omega$ is invertible, (2) $\Omega$ maps orbits to orbits and (3) $\Omega_*\nu_1 = \nu_2$. These are proven in the next three sections.



 \section{The inverse of $\Omega$}\label{sec:invert}
 In this section, we prove:
 \begin{prop}\label{prop:inverse}
 The map $\Omega:\cG_1 \to \cG_2$ is invertible with measurable inverse.
 \end{prop}
To prove this we will explicitly construct the inverse by swapping $\Gamma_1$ with $\Gamma_2$ in the construction of $\Omega$. Let $\upsilon:G_2 \times X \to G_1$ be the Zimmer 1-cocycle satisfying
\begin{eqnarray}\upsilon(g,x)\cdot x = g\cdot x,~\forall g\in G_2, x\in X,\label{eqn:cocycle3}\\
\upsilon(g_1g_2,x) = \upsilon(g_1, g_2\cdot x)\upsilon(g_2,x),~\forall g_1,g_2\in G_2, x\in X.\label{eqn:cocycle4}\end{eqnarray}
Note that
\begin{eqnarray}\label{eqn:cocycle5}
\omega(\upsilon(g_2,x),x) = g_2, ~ \upsilon(\omega(g_1,x),x)=g_1,~ \forall g_1 \in G_1, g_2\in G_2, x\in X.
\end{eqnarray}

 \begin{lem}
The range of $\Omega$ is contained in $\cG_2$. 
\end{lem}

\begin{proof}
Let $f\in \cG_1$, $g\in G_2$ and $\gamma \in \Gamma_1$. It suffices to show that
$$(\Omega f)\big(\beta(\gamma^{-1},f)^{-1}g_2\big)  = g_2^{-1}\cdot (\Omega f)\big(\beta(\gamma^{-1},f\big)^{-1}).$$
Observe:
\begin{eqnarray*}
g_2^{-1} \cdot (\Omega f)\big(\beta(\gamma^{-1},f)^{-1}\big) &=& g_2^{-1} \cdot f(\gamma)\\
&=& \upsilon(g_2^{-1}, f(\gamma)) \cdot f(\gamma)\\
&=& f(\gamma \upsilon(g_2^{-1}, f(\gamma))^{-1})\\
&=& (\Omega f)\Big(\beta(  \upsilon\big(g_2^{-1}, f(\gamma)\big)\gamma^{-1}, f)^{-1}\Big).
\end{eqnarray*}
So it suffices to show that
$$\beta(\gamma^{-1},f)^{-1}g_2 = \beta\big(  \upsilon(g_2^{-1}, f(\gamma)\big)\gamma^{-1}, f)^{-1}.$$
By taking inverses of both sides we see that this is equivalent to
$$g_2^{-1}\beta(\gamma^{-1},f) = \beta\Big(  \upsilon\big(g_2^{-1}, f(\gamma)\big)\gamma^{-1}, f\Big).$$
Proposition \ref{prop:beta} implies
$$\beta\Big(  \upsilon\big(g_2^{-1}, f(\gamma)\big)\gamma^{-1}, f\Big) = \omega\Big( \upsilon\big(g_2^{-1}, f(\gamma)\big), f(\gamma)\Big) \beta( \gamma^{-1}, f).$$
So it suffices to show that
$$g_2^{-1} = \omega\big( \upsilon(g_2^{-1}, f(\gamma)), f(\gamma)\big)$$
which is true by (\ref{eqn:cocycle5}).
\end{proof}


The proof of the next proposition is similar to that of proposition \ref{prop:beta}.
\begin{prop}\label{prop:alpha}
There exists a unique Borel map $\alpha:\Gamma_2 \times \cG_2 \to \Gamma_1$ satisfying the following.
\begin{itemize}
\item $\alpha(g,f) = \upsilon(g,f(e)),~ \forall g\in G_2$.
\item $\alpha(h,f) = h,~\forall h\in H$.
\item $\alpha(h\gamma, f) = h\alpha(\gamma,f), ~\forall h\in H, \gamma \in \Gamma_2$.
\item $\alpha(g\gamma,f) = \alpha(g,\gamma f) \alpha(\gamma,f) = \upsilon(g, f(\gamma^{-1}))\alpha(\gamma,f), ~\forall g\in G_2,\gamma \in \Gamma_2$.
\end{itemize}
\end{prop}


As in lemma \ref{lem:bijective}, $\alpha(\cdot,f)$ is bijective for any $f\in \cG_2$. So we define $\Theta:\cG_2 \to \cG_1$ by 
$$(\Theta f)(\alpha(\gamma^{-1},f)^{-1}) = f(\gamma), ~\forall f\in\cG_2, \gamma \in \Gamma_2.$$

We will show that $\Theta$ is the inverse of $\Omega$. We will need the next lemma.

\begin{lem}\label{lem:inverse}
The following are true:
\begin{eqnarray*}
\alpha\big(\beta(\gamma,f), \Omega f\big) &=& \gamma, ~\forall \gamma \in \Gamma_1, f\in \cG_1\\
\beta\big(\alpha(\gamma,f), \Theta f\big) &=& \gamma, ~\forall \gamma \in \Gamma_2, f\in \cG_2.
\end{eqnarray*}
\end{lem}

\begin{proof}
We prove only the first equation since the second one is similar. It is easy to check that the lemma is true if $\gamma=e$. By induction, it suffices to show that for any $\gamma \in\Gamma_1$, $g\in G_1$ and $h\in H$,
\begin{eqnarray*}
\alpha\big(\beta(g\gamma,f), \Omega f\big)&=&g\alpha\big(\beta(\gamma,f), \Omega f\big);\\
\alpha\big(\beta(h\gamma,f), \Omega f\big)&=&h\alpha\big(\beta(\gamma,f), \Omega f\big)
\end{eqnarray*}
The second equation above is immediate. The first equation is a short calculation:
\begin{eqnarray*}
\alpha\big(\beta(g\gamma,f), \Omega f\big)&=& \alpha\Big(\omega\big(g,f(\gamma^{-1})\big)\beta(\gamma,f),\Omega f\Big)\\
&=& \upsilon\Big( \omega\big(g,f(\gamma^{-1})\big), \Omega f\big( \beta(\gamma,f)^{-1}\big)\Big)\alpha\big(\beta(\gamma,f), \Omega f\big)\\
&=& \upsilon\Big( \omega\big(g,f(\gamma^{-1})\big), f(\gamma^{-1})\Big)\alpha\big(\beta(\gamma,f), \Omega f\big)\\
&=&g \alpha\big(\beta(\gamma,f), \Omega f\big).
\end{eqnarray*}
The first equality uses proposition \ref{prop:beta}, the second uses proposition \ref{prop:alpha}, the third uses the definition of $\Omega$ and the last uses equation (\ref{eqn:cocycle5}). 
\end{proof}

\begin{lem}\label{lem:weird}
The following are true:
\begin{eqnarray*}
\Omega f(\gamma) &=& f(\alpha(\gamma^{-1},\Omega f)^{-1}), ~\forall f\in \cG_1, \gamma \in \Gamma_2\\
\Theta f(\gamma) &=& f(\beta(\gamma^{-1},\Theta f)^{-1}), ~\forall f\in \cG_2, \gamma \in \Gamma_1.
\end{eqnarray*}
\end{lem}

\begin{proof}
We prove only the first equation since the second one is similar. By definition, if $\tau \in \Gamma_1$ is such that $\gamma= \beta(\tau^{-1},f)^{-1}$ then $\Omega f(\gamma) = f(\tau).$ The previous lemma implies that
$$\alpha\big(\beta(\tau^{-1},f),\Omega f\big) = \tau^{-1}.$$
Therefore, $\tau^{-1} = \alpha(\gamma^{-1},\Omega f)$ and 
$$\Omega f(\gamma) = f(\tau) = f\big( \alpha(\gamma^{-1},\Omega f)^{-1} \big)$$
as claimed.
\end{proof}

\begin{lem}\label{lem:thetaomega}
For any $f\in \cG_1$, $\Theta\Omega(f)=f$. Also, for any $f\in \cG_2$, $\Omega\Theta(f)=f$.
\end{lem}

\begin{proof}
For any $f \in \cG_1$ and any $\gamma \in \Gamma_2$,
$$\Theta \Omega f\big(\alpha(\gamma^{-1}, \Omega f)^{-1}\big) = \Omega f(\gamma) = f(\alpha(\gamma^{-1},\Omega f)^{-1}).$$
Therefore, $\Theta \Omega f = f$ as claimed. The second statement can be proven similarly.
\end{proof}
Proposition \ref{prop:inverse} follows immediately from the lemma above.

\section{Orbits to orbits}\label{sec:orbits}

\begin{lem}\label{lem:cocycle}
The cocycles $\beta$ and $\alpha$ satisfy the equations:
\begin{eqnarray*}
\beta(\gamma_1\gamma_2,f) &=& \beta(\gamma_1,\gamma_2 f) \beta(\gamma_2,f),~\forall \gamma_1,\gamma_2\in \Gamma_1, f\in \cG_1\\
\alpha(\gamma_1\gamma_2,f) &=& \alpha(\gamma_1,\gamma_2 f) \alpha(\gamma_2,f),~\forall \gamma_1,\gamma_2\in \Gamma_2, f\in \cG_2.
\end{eqnarray*}
\end{lem}
\begin{proof}
This follows from propositions \ref{prop:beta} and \ref{prop:alpha}. 
\end{proof}

\begin{lem}
For any $\gamma \in \Gamma_1, f\in \cG_1$,
$$\beta(\gamma,f)(\Omega f) = \Omega(\gamma f).$$
Similarly,
$$\alpha(\gamma,f)(\Theta f) = \Theta(\gamma f),~\forall \gamma \in \Gamma_2, f\in \cG_2.$$
\end{lem}
\begin{proof}
We prove only the first equation since the second one is similar. By lemma \ref{lem:bijective}, it suffices to show that
$$\beta(\gamma,f)(\Omega f)\big(\beta(\tau^{-1},\gamma f)^{-1}\big) = \Omega(\gamma f)\big(\beta(\tau^{-1},\gamma f)^{-1}\big)~\forall \tau \in \Gamma_1.$$
By definition of the action of $\Gamma_2$ on $\cG_2$:
$$ \beta(\gamma,f)(\Omega f)\big(\beta(\tau^{-1},\gamma f)^{-1}\big)  = (\Omega f)\big(\beta(\gamma,f)^{-1}\beta(\tau^{-1},\gamma f)^{-1}\big).$$
By the previous lemma,
$$\beta(\gamma,f)^{-1}\beta(\tau^{-1},\gamma f)^{-1} = \beta(\tau^{-1}\gamma,f)^{-1}.$$
By definition of $\Omega f$
$$(\Omega f)\big(\beta(\tau^{-1}\gamma,f)^{-1}\big) = f(\gamma^{-1}\tau) = (\gamma f)(\tau) =\Omega(\gamma f)\big(\beta(\tau^{-1},\gamma f)^{-1}\big).$$
\end{proof}

The lemma above and lemma \ref{lem:bijective} imply:
\begin{prop}\label{prop:orbits}
For every $f\in \cG_1$,  $\Omega \big(\{\gamma f:~ \gamma \in \Gamma_1\big\}\big) = \{\gamma (\Omega f):~\gamma \in \Gamma_2\}.$
\end{prop}

\section{Measure-space isomorphism}\label{sec:measure}

\begin{prop}\label{prop:measure}
The pushforward measure $\Omega_*\nu_1 = \nu_2$.
\end{prop}


Let us recall the definition of the measures $\nu_i$ for $i=1,2$. Let $\sigma_i:\Gamma_i/G_i \to \Gamma_i$ be a section such that $\sigma_i(G_i)=e$. Define $\Phi_i:\cG_i\to X^{\Gamma_i/G_i}$ by $\Phi_i(f)(C)=f\big(\sigma_i(C)\big)$. $\Phi_i$ is a bijection. By definition, $\nu_i$ is the pullback measure $\Phi_i^*\mu^{\Gamma_i/G_i}$. 

The measures $\nu_i$ do not depend on the choice of section. So, we will make the following choice. For each coset $C\in \Gamma_i/G_i$ such that $C\ne G_i$, let $\sigma_i(C)=g_1h_1\cdots g_{n}h_n$ where $g_1h_1\cdots g_nh_n$ is the unique element of $C$ such that $g_1 \in G_1$, $g_2,\ldots, g_n \in G_1-\{e\}$ and $h_1,\ldots, h_n \in H-\{e\}$. Define length$(C)=n$. Also let length$(G_i)=0$ and if $\gamma \in \Gamma_i$ define length$(\gamma):=$length$(\gamma G_i)$. 

Let $\rng \sigma_i :=\sigma_i(\Gamma_i/G_i)$. Let $f \in \cG_1$ be chosen at random with law $\nu_1$. For $\gamma \in \Gamma_2$, let $Z_\gamma :=\Omega f(\gamma)$. Proposition \ref{prop:measure} is equivalent to the assertion that $\{Z_\gamma:~ \gamma \in \rng \sigma_2\}$ is a jointly independent set of random variables and each $Z_\gamma \in X$ has distribution $\mu$ (this is because $\Omega f \in \cG_2$, so its values on $\Gamma_2$ are determined by its values on $\rng \sigma_2$). Note that $Z_e = \Omega f(e) = f(e)$ has distribution $\mu$. Because the pushforward measure $\Omega_*\nu_1$ is $\Gamma_2$-invariant, it follows that every $Z_\gamma$ has distribution $\mu$. So it suffices to show that the family $\{Z_\gamma:~\gamma \in \rng \sigma_2\}$ is jointly independent.

\begin{lem}\label{lem:basic}
For any $f\in \cG_1$, the map $B_f:\Gamma_1 \to \Gamma_2$ defined by $B_f(\gamma) :=\beta(\gamma^{-1},f)^{-1}$ restricts to a bijection from $\rng \sigma_1$ to $\rng \sigma_2$. Moreover, length$(B_f(\gamma))=$length$(\gamma)$. Similarly,  the map $A_{\Omega f}:\Gamma_2 \to \Gamma_1$ defined by $A_{\Omega f}(\gamma) :=\alpha(\gamma^{-1},\Omega f)^{-1}$ restricts to a bijection from $\rng \sigma_2$ to $\rng \sigma_1$ and length$(A_{\Omega f}(\gamma))=$length$(\gamma)$.
 \end{lem}

\begin{proof}
Proposition \ref{prop:beta} implies $B_f$ maps  $\rng \sigma_1$ into $\rng \sigma_2$ and preserves lengths. Similarly, proposition \ref{prop:alpha} implies the map $A_{\Omega f}$ maps $\rng \sigma_2$ into $\rng \sigma_1$ and preserves lengths. Lemmas \ref{lem:inverse} and \ref{lem:thetaomega} imply these maps are inverses.
\end{proof}

For $\gamma \in \Gamma_1$, let $Y_\gamma := f(\gamma)$ (where, as above, $f \in \cG_1$ is chosen at random with law $\nu_1$). By definition of $\nu_1$, $\{Y_\gamma:~\gamma \in \rng \sigma_1\}$ is an i.i.d. (independent identically distributed) family. It is tempting then to claim that, since $Z_{\beta(\gamma^{-1},f)^{-1}} = Y_\gamma$,  the lemma above implies $\{Z_\gamma:~\gamma \in \rng \sigma_2\}$ is an i.i.d. family. However, since the map $\gamma \mapsto \beta(\gamma^{-1},f)^{-1}$ depends on $f$, this cannot be taken for granted. 

Let $L_i(n) = \{\sigma_i(C) \in \Gamma_i:~ C \in \Gamma_i/G_i, ~$length$(C) \le n\}$. Assume for induction that for some $n\ge 0$, $\{Z_\gamma:~\gamma \in L_2(n)\}$ is i.i.d. (this is trivially true if $n=0$). By lemma \ref{lem:weird}, $Z_\gamma = Y_{\alpha(\gamma^{-1}, \Omega f)^{-1}}$. The previous lemma implies that 
$$\{ Y_{\alpha(\gamma^{-1}, \Omega f)^{-1}}:~\gamma \in L_2(n+1)\} = \{Y_\gamma:~ \gamma \in L_1(n+1)\}$$
is i.i.d.. So it suffices to show that the map 
$$\gamma \in L_2(n+1) - L_2(n),~ \gamma \mapsto \alpha(\gamma^{-1}, \Omega f)^{-1}$$
depends only on the variables $\{Y_\gamma:~ \gamma \in L_1(n)\}$ (i.e., that it is a function of these variables only). This is accomplished in the next lemma.

\begin{lem}
Let $n \ge 0$. Let $f_1,f_2 \in \cG_1$ and suppose that $f_1(\gamma)=f_2(\gamma)$ for all $\gamma \in L_1(n)$. Then
$$\alpha(\gamma^{-1},\Omega f_1) = \alpha(\gamma^{-1},\Omega f_2),~\forall \gamma \in L_2(n+1).$$
\end{lem}
\begin{proof}
If $\gamma \in L_2(0)$ then $\gamma=e$ and the result is trivial. So assume for induction that the equation holds for all $\gamma \in L_2(k)$ for some $0\le k \le n$. Let $\gamma = g_1h_1\cdots g_{k+1}h_{k+1} \in L_2(k+1)$. Let $\tau = g_1h_1\cdots g_{k}h_{k} \in L_2(k)$. By lemma \ref{lem:cocycle} and proposition \ref{prop:alpha},
\begin{eqnarray*}
\alpha(\gamma^{-1},\Omega f_1) &=&  \alpha(h_{k+1}^{-1}g_{k+1}^{-1}\tau^{-1},\Omega f_1)\\
&=& \alpha(h_{k+1}^{-1}g_{k+1}^{-1},\tau^{-1}\Omega f_1)  \alpha(\tau^{-1},\Omega f_1)\\
&=&h_{k+1}^{-1}\upsilon\big(g_{k+1}^{-1}, \tau^{-1}\Omega f_1(e)\big) \alpha(\tau^{-1},\Omega f_1).
\end{eqnarray*}
Similar reasoning shows that 
$$\alpha(\gamma^{-1},\Omega f_2) =h_{k+1}^{-1}\upsilon\big(g_{k+1}^{-1}, \tau^{-1}\Omega f_2(e)\big) \alpha(\tau^{-1},\Omega f_2).$$
The induction hypothesis implies $\alpha(\tau^{-1},\Omega f_1) = \alpha(\tau^{-1},\Omega f_2)$. So it suffices to show that $ \tau^{-1}\Omega f_1(e) = \tau^{-1}\Omega f_2(e)$. By lemma \ref{lem:weird},
$$\tau^{-1}\Omega f_1(e) = \Omega f_1(\tau) = f_1(\alpha(\tau^{-1},\Omega f_1)^{-1})= f_1(\alpha(\tau^{-1},\Omega f_2)^{-1}).$$
By the previous lemma, $\alpha(\tau^{-1},\Omega f_2)^{-1} \in L_1(k) \subset L_1(n)$. So the hypotheses on $f_1$ and $f_2$ imply
$$\tau^{-1}\Omega f_1(e)=f_1(\alpha(\tau^{-1},\Omega f_2)^{-1}) = f_2(\alpha(\tau^{-1},\Omega f_2)^{-1})=\tau^{-1}\Omega f_2(e).$$

\end{proof}
The lemma above and the preceding discussion show that $\{Z_\gamma:~\gamma \in L_2(n)\}$ is i.i.d. for all $n\ge 0$. This implies that the family $\{Z_\gamma:~\gamma \in \Gamma_2\}$ is i.i.d. which implies proposition \ref{prop:measure}. Theorem \ref{thm:main} follows immediately from propositions \ref{prop:inverse}, \ref{prop:orbits} and \ref{prop:measure}.


\end{document}